# A special class of Hankel determinants

Johann Cigler

Fakultät für Mathematik, Universität Wien
johann.cigler@univie.ac.at

**Abstract**

In this expository paper we compute Hankel determinants of some sequences whose generating functions are given by C- fractions and derive orthogonality properties for associated polynomials.

## 1. Introduction

Let $(b_k)_{k\geq -1}$ be a non-decreasing sequence of integers such that $b_{k+2}-b_k \geq 1$ with initial values $b_{-1}=-1$ and $b_0=0$ and define the formal power series $f(x)$ by the C-fraction

$$f(x)=\sum_{n\geq 0} f_n x^n = \cfrac{1}{1-\cfrac{a_0 x^{b_1-b_{-1}}}{1-\cfrac{a_1 x^{b_2-b_0}}{1-\cfrac{a_2 x^{b_3-b_1}}{1-\cfrac{a_3 x^{b_4-b_2}}{1-\ddots}}}}}. \tag{1.1}$$

The Hankel determinants of the coefficients $f_n$ will be denoted by $d(n)=\det\left(f_{i+j}\right)_{i,j=0}^{n}$.
Then all non-vanishing Hankel determinants are given by

$$d(b_k)=(-1)^{\binom{b_1-b_0}{2}+\binom{b_2-b_1}{2}+\cdots\binom{b_k-b_{k-1}}{2}} a_0^{b_k-b_0} a_1^{b_k-b_1} a_2^{b_k-b_2} a_3^{b_k-b_3}\cdots a_{k-1}^{b_k-b_{k-1}}. \tag{1.2}$$

For $b_k=k$ or $b_k=\left\lfloor \frac{k}{2}\right\rfloor$, i.e. $b_k-b_{k-2}=2$ or $b_k-b_{k-2}=1$ these are old results which are intimately connected with orthogonal polynomials. For the more general case Paul Barry [2] conjectured that $d(b_k)=\pm 1$ if all $a_n=\pm 1.$ This conjecture led me to investigate this problem. Afterwards I learned that an equivalent result had already earlier been found by Victor I. Buslaev [4] with another proof.

If $b_{k+1}-b_k>1$ for some $k$ the role of orthogonal polynomials is played by polynomials $r_k(x)$ which satisfy the three-term recurrence $r_k(x)=x^{b_{k-1}-b_{k-2}}r_{k-1}(x)-a_{k-2}r_{k-2}(x)$ with initial values $r_0(x)=1$ and $r_1(x)=x.$ They satisfy $\Lambda\left(r_k(x)x^n\right)=0$ for $n<b_k$ and $\Lambda\left(r_k(x)x^{b_k}\right)=a_0\cdots a_{k-1},$ if $\Lambda$ denotes the linear functional on the vector space of polynomials defined by $\Lambda\left(x^n\right)=f_n.$
We show that the polynomials $r_k(x)$ are related to the polynomials $p_n(x)=\det\left(f_{i+j}x-f_{i+j+1}\right)_{i,j=0}^{n-1}$ and study some typical examples.

## 2. Some preliminary results about continued fractions and Hankel determinants

Let us begin with some well-known facts about continued fractions.

If we write

$$\frac{1}{c_0+\cfrac{a_0}{c_1+\cfrac{a_1}{c_2+\ddots\cfrac{}{c_{n-1}}}}}=\frac{1}{c_0+}\frac{a_0}{c_1+}\cdots\frac{a_{n-2}}{c_{n-1}}=\frac{A_n}{B_n} \tag{2.1}$$

then $A_0=0, A_1=1$ and $B_0=1, B_1=c_0$
and

$$\begin{aligned} A_n &= c_{n-1}A_{n-1}+a_{n-2}A_{n-2}, \\ B_n &= c_{n-1}B_{n-1}+a_{n-2}B_{n-2}. \end{aligned} \tag{2.2}$$

Moreover we get by induction

$$A_{k+1}B_k-B_{k+1}A_k=(-1)^k a_0a_1\cdots a_{k-1}. \tag{2.3}$$

Let us assume that $b_{-1}=-1,$ $b_0=0,$ and that $b_n$ is not decreasing with $b_{n+2}-b_n\geq 1.$ If $b_n=b_{n+1}=b$ for some $n$ we say that $b$ occurs with multiplicity $2.$ Otherwise it has multiplicity $1.$ Since $b_{n+2}-b_n\geq 1$ there are no triplets $b_n=b_{n+1}=b_{n+2}.$
If each non-negative integer occurs in the sequence $(b_k)$ the above results are intimately connected with orthogonal polynomials. This fact has been generalized by E. Frank [7] to more general cases.

The above choice of $(b_k)$ leads to the following facts:
Let

$$f_n(x)=\frac{1}{1-\cfrac{a_0x^{b_1-b_{-1}}}{1-\cfrac{a_1x^{b_2-b_0}}{1-\ddots\cfrac{}{1-a_{n-2}x^{b_{n-1}-b_{n-3}}}}}}=\frac{1}{1-}\frac{a_0x^{b_1-b_{-1}}}{1-}\frac{a_1x^{b_2-b_0}}{1-}\cdots\frac{a_{n-3}x^{b_{n-2}-b_{n-4}}}{1-}\frac{a_{n-2}x^{b_{n-1}-b_{n-3}}}{1}=\frac{A_n(x)}{B_n(x)} \tag{2.4}$$

with $f_0(x)=0$ and $f_1(x)=1.$

Here (2.2) reduces to

$$A_n(x)=A_{n-1}(x)-a_{n-2}x^{b_{n-1}-b_{n-3}}A_{n-2}(x) \tag{2.5}$$

and

$$B_n(x)=B_{n-1}(x)-a_{n-2}x^{b_{n-1}-b_{n-3}}B_{n-2}(x) \tag{2.6}$$

with initial values $A_0(x)=0$, $A_1(x)=1$ and $B_0(x)=1, B_1(x)=1.$

As a generalization of orthogonal polynomials we define polynomials $r_k(x)$ by

$$r_k(x) = x^{b_{k-1}-b_{k-2}} r_{k-1}(x) - a_{k-2} r_{k-2}(x) \tag{2.7}$$

with initial values $r_0(x) = 1$ and $r_1(x) = x.$

It is easily verified that $\deg r_k(x) = b_{k-1} + 1$ for $k \geq 0.$

**Proposition 2.1**

*The polynomials* $r_k(x)$ *and* $B_k(x)$ *are related by*

$$r_k(x) = x^{b_{k-1}+1} B_k\left(\frac{1}{x}\right). \tag{2.8}$$

**Proof**

Let $U_k(x) = x^{b_{k-1}+1} B_k\left(\frac{1}{x}\right).$ By (2.6) we get

$U_k(x) = x^{b_{k-1}+1} B_k\left(\frac{1}{x}\right) = x^{b_{k-1}+1} B_{k-1}\left(\frac{1}{x}\right) - a_{k-2} x^{b_{k-1}+1} x^{-b_{k-1}+b_{k-3}} B_{k-2}\left(\frac{1}{x}\right)$

$= x^{b_{k-1}-b_{k-2}} U_{k-1}(x) - a_{k-2} U_{k-2}(x)$

with initial values $U_0(x) = 1$ and $U_1(x) = x.$

Therefore $U_k(x) = r_k(x).$

It is clear that $\deg A_n(x) \leq \deg B_n(x)$. By induction we get

$\deg B_n(x) \leq b_{n-1} + 1.$

More precisely we get $\deg B_{2n} = 1 + b_{2n-1}.$ For

$\deg\left(B_{2n}(x)\right) = \max\left(\deg B_{2n-1}(x), \deg x^{b_{2n-1}-b_{2n-3}} B_{2n-2}(x)\right) = 1 + b_{2n-1}.$

If for some $m$ $b_{2m-1} = b_{2m}$ then $\deg B_{2m+1}(x) = 1 + b_{2m-1} = 1 + b_{2m}$ and for all $n \geq m$ also $\deg B_{2n+1}(x) = 1 + b_{2n}.$ For $n < m$ we have $\deg B_{2n+1}(x) = b_{2n}.$ In this case $r_{2n+1}(0) = 0.$

This follows from

$B_{2n+1}(x) = B_{2n}(x) - a_{2n-1} x^{b_{2n}-b_{2n-2}} B_{2n-1}(x)$

by induction.

**Proposition 2.2**

*Let* $f(x) = \sum_{n \geq 0} f_n x^n$ *satisfy (3.1) and let*

$$B_k(x) f(x) - A_k(x) = \sum_k f_n^{(k)} x^n. \tag{2.9}$$

*Then* $f_n^{(k)} = 0$ *for* $n \leq b_{k-1} + b_k$ *and* $f_{b_{k-1}+b_k+1}^{(k)} = a_0 a_1 \cdots a_{k-1}.$

**Proof**

Observe that by (2.3)

$$A_{k+1}(x)B_k(x)-B_{k+1}(x)A_k(x)=a_0a_1\cdots a_{k-1}x^{1+b_{k-1}+b_k}. \tag{2.10}$$

Therefore we get the series expansion

$$\frac{A_{k+1}(x)}{B_{k+1}(x)}-\frac{A_k(x)}{B_k(x)}=\frac{a_0a_1\cdots a_{k-1}x^{1+b_{k-1}+b_k}}{B_{k+1}(x)B_k(x)}=a_0a_1\cdots a_{k-1}x^{1+b_{k-1}+b_k}+\cdots$$

since $B_k(0)=1$ for all $k$.
Since $1+b_{k-1}+b_k>1+b_{k-2}+b_{k-1}$ the series expansion of $f(x)$ coincides with the expansion of $\frac{A_k(x)}{B_k(x)}$ for the first $b_{k-1}+b_k$ terms.
Therefore the coefficients of $x^n$ in $f(x)-\frac{A_k(x)}{B_k(x)}$ vanish for
$n<b_1-b_{-1}+b_2-b_0+\cdots+b_k-b_{k-2}=b_{k-1}+b_k+1$.

Since $\deg A_k(x)\le b_{k-1}$ we see that the coefficients of $B_k(x)f(x)$ vanish for $b_{k-1}<n\le b_{k-1}+b_k$.

If we set $B_k(x)=\sum_j v_{k,j}x^j$. Then we get

$$v_{k,0}f_n+v_{k,1}f_{n-1}+\cdots=0 \text{ for } b_{k-1}<n\le b_{k-1}+b_k$$

and

$$v_{k,0}f_{b_{k-1}+b_k+1}+v_{k,1}f_{b_{k-1}+b_k}+\cdots=a_0a_1\cdots a_{k-1}.$$

Another equivalent formulation gives

**Proposition 2.3**

*Let $\Lambda$ denote the linear functional on the polynomials defined by*

$$\Lambda\left(x^n\right)=f_n. \tag{2.11}$$

*Then*

*for* $n<b_k$

$$\Lambda\left(r_k(x)x^n\right)=\Lambda\left(x^{b_{k-1}+1}B_k\left(\frac{1}{x}\right)x^n\right)=\Lambda\left(\sum_j v_{k,j}x^{b_{k-1}-j+n+1}\right)=\sum_j v_{k,j}f_{b_{k-1}-j+n+1}=0 \tag{2.12}$$

*and*

$$\Lambda\left(r_k(x)x^{b_k}\right)=\Lambda\left(\sum_j v_{k,j}x^{b_{k-1}-j+b_k+1}\right)=\sum_j v_{k,j}f_{b_{k-1}-j+b_k+1}=a_0\cdots a_{k-1}. \tag{2.13}$$

**Remark**

The identities (2.12) and (2.13) generalize the concept of orthogonality and (2.7) is an analogue of the three-term recurrence of orthogonal polynomials.

For some polynomials $r_k(x)$ a simple formula can be found. To this end we define the polynomials

$$p_n(x) = \det\begin{pmatrix} f_0 & f_1 & \cdots & f_{n-1} & 1 \\ f_1 & f_2 & \cdots & f_n & x \\ f_2 & f_3 & \cdots & f_{n+1} & x^2 \\ \vdots & & & & \vdots \\ f_n & f_{n+1} & \cdots & f_{2n-1} & x^n \end{pmatrix} \tag{2.14}$$

for $n \geq 1$ and let $p_0(x) = 1$.

Note that by elementary row operations $p_n(x)$ can also be expressed as

$$p_n(x) = \det\left( f_{i+j}x - f_{i+j+1}\right)_{i,j=0}^{n-1}. \tag{2.15}$$

We will need some generalizations of

**Proposition 2.4 (G.E.Andrews and J. Wimp [1])**

*Let*

$$s(x) = \sum_{n\geq 0} s_n x^n \tag{2.16}$$

*with* $s_0 = 1$ *and*

$$t(x) = \frac{1}{s(x)} = \sum_{n\geq 0} t_n x^n. \tag{2.17}$$

*Then for* $n \geq 1$

$$\det\left(s_{i+j}\right)_{i,j=0}^{n} = (-1)^n \det\left(t_{i+j+2}\right)_{i,j=0}^{n-1}. \tag{2.18}$$

In view of our applications we note the following special case:

**Corollary 2.1**

$$f(x) = \sum_{n\geq 0} f_n x^n = \frac{1}{1 - ax^2 \sum_{j\geq 0} c_j x^j} \tag{2.19}$$

*implies*

$$\det\left( f_{i+j}\right)_{i,j=0}^{n} = a^n \det\left(c_{i+j}\right)_{i,j=0}^{n-1}. \tag{2.20}$$

**Example 2.1**

*Let*

$$f(x,a_0,a_1,a_2,\cdots)=\sum_{n\geq 0} f_n x^n = \cfrac{1}{1-\cfrac{a_0x^2}{1-\cfrac{a_1x^2}{1-\cfrac{a_2x^2}{1-\cfrac{a_3x^2}{1-\ddots}}}}}. \tag{2.21}$$

*Then*

$$d(n)=a_0^n a_1^{n-1} a_2^{n-2}\cdots a_{n-1}. \tag{2.22}$$

This well-known result follows from $f(x,a_0,a_1,a_2,\cdots)=\dfrac{1}{1-a_0x^2 f(x,a_1,a_2,\cdots)}.$

**Example 2.1.1**

If $a_n=1$ for all $n$ then $f(x)$ satisfies $1-f(x)+x^2f(x)^2=0$ which gives

$$f(x)=\frac{1-\sqrt{1-4x^2}}{2x^2}=\sum_{n\geq 0} C_n x^{2n} \tag{2.23}$$

with the Catalan numbers

$$f_{2n}=C_n=\frac{1}{n+1}\binom{2n}{n}. \tag{2.24}$$

**Proposition 2.5**

*Let* $s(x)=\sum_{n\geq 0} s_n x^n$ *with* $s_0=1$ *and* $t(x)=\dfrac{1}{s(x)}=\sum_{n\geq 0} t_n x^n.$

*Define* $s_n=0$ *and* $\det\left(t_{i+j+k}\right)_{i,j=0}^{n}=1$ *for* $n<0.$ *Let* $m\geq 0.$

*Then*

$\det\left(s_{i+j-m}\right)_{i,j=0}^{n}=0$ *for* $n<m$ *and*

$$\det\left(s_{i+j-m}\right)_{i,j=0}^{n+m}=(-1)^{n+\binom{m+1}{2}}\det\left(t_{i+j+m+2}\right)_{i,j=0}^{n-1}. \tag{2.25}$$

The proof is almost the same as in [1]. I shall illustrate the method (called O-reduction in [1]) for $m=2$ and $n=3.$ Observe that $t_0=1$ and $\sum_{j=0}^{n} s_{n-j}t_j=[n=0].$

$$\det\left(s_{i+j-2}\right)_{i,j=0}^{5}=\det\begin{pmatrix} 0 & 0 & s_0 & s_1 & s_2 & s_3 \\ 0 & s_0 & s_1 & s_2 & s_3 & s_4 \\ s_0 & s_1 & s_2 & s_3 & s_4 & s_5 \\ s_1 & s_2 & s_3 & s_4 & s_5 & s_6 \\ s_2 & s_3 & s_4 & s_5 & s_6 & s_7 \\ s_3 & s_4 & s_5 & s_6 & s_7 & s_8 \end{pmatrix}$$

Then we get with elementary column operations

$$\det\left(s_{i+j-2}\right)_{i,j=0}^{5}=\det\begin{pmatrix} 0 & 0 & s_0 & t_0s_1+t_1s_0 & t_0s_2+t_1s_1+t_2s_0 & t_0s_3+t_1s_2+t_2s_1+t_3s_0 \\ 0 & s_0 & s_1 & t_0s_2+t_1s_1 & t_0s_3+t_1s_2+t_2s_1 & t_0s_4+t_1s_3+t_2s_2+t_3s_1 \\ s_0 & s_1 & s_2 & t_0s_3+t_1s_2 & t_0s_4+t_1s_3+t_2s_2 & t_0s_5+t_1s_4+t_2s_3+t_3s_2 \\ s_1 & s_2 & s_3 & t_0s_4+t_1s_3 & t_0s_5+t_1s_4+t_2s_3 & t_0s_6+t_1s_5+t_2s_4+t_3s_3 \\ s_2 & s_3 & s_4 & t_0s_5+t_1s_4 & t_0s_6+t_1s_5+t_2s_4 & t_0s_7+t_1s_6+t_2s_5+t_3s_4 \\ s_3 & s_4 & s_5 & t_0s_6+t_1s_5 & t_0s_7+t_1s_6+t_2s_5 & t_0s_8+t_1s_7+t_2s_6+t_3s_5 \end{pmatrix}$$

$$=\det\begin{pmatrix} 0 & 0 & s_0 & 0 & 0 & 0 \\ 0 & s_0 & s_1 & -t_2s_0 & -t_3s_0 & -t_4s_0 \\ s_0 & s_1 & s_2 & -t_2s_1-t_3s_0 & -t_3s_1-t_4s_0 & -t_4s_1-t_5s_0 \\ s_1 & s_2 & s_3 & -t_2s_2-t_3s_1-t_4s_0 & -t_3s_2-t_4s_1-t_5s_0 & -t_4s_2-t_5s_1-t_6s_0 \\ s_2 & s_3 & s_4 & -t_2s_3-t_3s_2-t_4s_1-t_5s_0 & -t_3s_3-t_4s_2-t_5s_1-t_6s_0 & -t_4s_3-t_5s_2-t_6s_1-t_7s_0 \\ s_3 & s_4 & s_5 & -t_2s_4-t_3s_3-t_4s_2-t_5s_1-t_6s_0 & -t_3s_4-t_4s_3-t_5s_2-t_6s_1-t_7s_0 & -t_4s_4-t_5s_3-t_6s_2-t_7s_1-t_8s_0 \end{pmatrix}$$

Finally with elementary row operations we get

$$=\det\begin{pmatrix} 0 & 0 & s_0 & 0 & 0 & 0 \\ 0 & s_0 & s_1 & -t_2s_0 & -t_3s_0 & -t_4s_0 \\ s_0 & 0 & s_2-\frac{s_1^2}{s_0} & -t_3s_0 & -t_4s_0 & -t_5s_0 \\ s_1 & 0 & s_3-\frac{s_2s_1}{s_0} & -t_3s_1-t_4s_0 & -t_4s_1-t_5s_0 & -t_5s_1-t_6s_0 \\ s_2 & 0 & s_4-\frac{s_3s_1}{s_0} & -t_3s_2-t_4s_1-t_5s_0 & -t_4s_2-t_5s_1-t_6s_0 & -t_5s_2-t_6s_1-t_7s_0 \\ s_3 & 0 & s_5-\frac{s_4s_1}{s_0} & -t_3s_3-t_4s_2-t_5s_1-t_6s_0 & -t_4s_3-t_5s_2-t_6s_1-t_7s_0 & -t_5s_3-t_6s_2-t_7s_1-t_8s_0 \end{pmatrix}$$

By iterating this prodecure we get finally

$$= \det \begin{pmatrix} 0 & 0 & s_0 & 0 & 0 & 0 \\ 0 & s_0 & s_1 & -t_2 s_0 & -t_3 s_0 & -t_4 s_0 \\ s_0 & 0 & * & -t_3 s_0 & -t_4 s_0 & -t_5 s_0 \\ 0 & 0 & * & -t_4 s_0 & -t_5 s_0 & -t_6 s_0 \\ 0 & 0 & * & -t_5 s_0 & -t_6 s_0 & -t_7 s_0 \\ 0 & 0 & * & -t_6 s_0 & -t_7 s_0 & -t_8 s_0 \end{pmatrix} = (-1)^{\binom{3}{2}+3} \det \begin{pmatrix} t_4 & t_5 & t_6 \\ t_5 & t_6 & t_7 \\ t_6 & t_7 & t_8 \end{pmatrix} = (-1)^{\binom{3}{2}+3} \det \left( t_{i+j+2+2} \right)_{i,j=0}^{3-1}.$$

An analogous result also holds for $m = -1.$

**Proposition 2.6**

*Let* $s(x) = \sum_{n \geq 0} s_n x^n$ *with* $s_0 = 1$ *and* $t(x) = \dfrac{1}{s(x)} = \sum_{n \geq 0} t_n x^n.$

*Then*

$$\det \left( s_{i+j+1} \right)_{i,j=0}^{n} = (-1)^{n+1} \det \left( t_{i+j+1} \right)_{i,j=0}^{n}. \tag{2.26}$$

**Proof**

I will illustrate the proof for $n = 2.$

$$\det \begin{pmatrix} s_1 & s_2 & s_3 \\ s_2 & s_3 & s_4 \\ s_3 & s_4 & s_5 \end{pmatrix} = \det \begin{pmatrix} s_1 & s_2 & -s_2 t_1 - s_1 t_2 - s_0 t_3 \\ s_2 & s_3 & -s_3 t_1 - s_2 t_2 - s_1 t_3 - s_0 t_4 \\ s_3 & s_4 & -s_4 t_1 - s_3 t_2 - s_2 t_3 - s_1 t_4 - s_0 t_5 \end{pmatrix}$$

$$= \det \begin{pmatrix} s_1 & s_2 & -s_0 t_3 \\ s_2 & s_3 & -s_1 t_3 - s_0 t_4 \\ s_3 & s_4 & -s_2 t_3 - s_1 t_4 - s_0 t_5 \end{pmatrix} = \det \begin{pmatrix} s_1 & -s_1 t_1 - s_0 t_2 & -s_0 t_3 \\ s_2 & -s_2 t_1 - s_1 t_2 - s_0 t_3 & -s_1 t_3 - s_0 t_4 \\ s_3 & -s_3 t_1 - s_2 t_2 - s_1 t_3 - s_0 t_4 & -s_2 t_3 - s_1 t_4 - s_0 t_5 \end{pmatrix}$$

$$= \det \begin{pmatrix} s_1 & -s_0 t_2 & -s_0 t_3 \\ s_2 & -s_1 t_2 - s_0 t_3 & -s_1 t_3 - s_0 t_4 \\ s_3 & -s_2 t_2 - s_1 t_3 - s_0 t_4 & -s_2 t_3 - s_1 t_4 - s_0 t_5 \end{pmatrix} = \det \begin{pmatrix} -s_0 t_1 & -s_0 t_2 & -s_0 t_3 \\ -s_1 t_1 - s_0 t_2 & -s_1 t_2 - s_0 t_3 & -s_1 t_3 - s_0 t_4 \\ -s_2 t_1 - s_1 t_2 - s_0 t_3 & -s_2 t_2 - s_1 t_3 - s_0 t_4 & -s_2 t_3 - s_1 t_4 - s_0 t_5 \end{pmatrix}.$$

The last determinant obviously equals

$$\det \begin{pmatrix} -s_0 t_1 & -s_0 t_2 & -s_0 t_3 \\ -s_0 t_2 & -s_0 t_3 & -s_0 t_4 \\ -s_0 t_3 & -s_0 t_4 & -s_0 t_5 \end{pmatrix} = -\det \begin{pmatrix} t_1 & t_2 & t_3 \\ t_2 & t_3 & t_4 \\ t_3 & t_4 & t_5 \end{pmatrix}.$$

**Corollary 2.2**

*Let*

$$f(x) = \frac{1}{1 - axc(x)}. \tag{2.27}$$

*Then*

$$\det\left( f_{i+j} \right)_{i,j=0}^{n} = a^n \det\left( c_{i+j+1} \right)_{i,j=0}^{n-1} \tag{2.28}$$

*and*

$$\det\left( f_{i+j+1} \right)_{i,j=0}^{n} = a^{n+1} \det\left( c_{i+j} \right)_{i,j=0}^{n}. \tag{2.29}$$

**Example 2.2**

$$f(x, a_0, a_1, a_2, \cdots) = \cfrac{1}{1 - \cfrac{a_0 x}{1 - \cfrac{a_1 x}{1 - \cfrac{a_2 x}{1 - \cdots}}}} \tag{2.30}$$

*implies*

$$\det\left( f_{i+j}\left(a_0, a_1, a_2, \cdots\right) \right)_{i,j=0}^{n} = \left(a_0 a_1\right)^n \left(a_2 a_3\right)^{n-1} \left(a_4 a_5\right)^{n-2} \cdots \left(a_{2n-2} a_{2n-1}\right). \tag{2.31}$$

**Proof**

Here we have

$$f(x, a_0, a_1, a_2, \cdots) = \frac{1}{1 - a_0 x f(x, a_1, a_2, \cdots)}$$

and

$$f(x, a_1, a_2, \cdots) = \frac{1}{1 - a_1 x f(x, a_2, a_3, \cdots)}.$$

Therefore

$$\det\left( f_{i+j}\left(a_0, a_1, a_2, \cdots\right) \right)_{i,j=0}^{n} = a_0^n \det\left( f_{i+j+1}\left(a_1, a_2, \cdots\right) \right)_{i,j=0}^{n-1}$$

and

$$\det\left( f_{i+j+1}\left(a_1, a_2, \cdots\right) \right)_{i,j=0}^{n-1} = a_1^n \det\left( f_{i+j}\left(a_2, a_3, \cdots\right) \right)_{i,j=0}^{n-1}.$$

This gives

$$\det\left( f_{i+j}\left(a_0, a_1, a_2, \cdots\right) \right)_{i,j=0}^{n} = \left(a_0 a_1\right)^n \det\left( f_{i+j}\left(a_2, a_3, a_4, \cdots\right) \right)_{i,j=0}^{n-1}.$$

By induction we get (2.31).

If all $a_n = 1$ then $f_n = C_n$.

These Propositions lead us to

**Lemma 2.1**

*Let* $m \geq -1$, $p \geq 1$ *and let* $f(x) = \sum_{n\geq 0} f_n x^n$ *and* $g(x) = \sum_{n\geq 0} g_n x^n$ *be formal power series with* $g_0 = 1$ *and define* $f_n = 0$ *for* $n < 0$.

*For*

$$f(x) = \frac{1}{1 - ax^p g(x)} \tag{2.32}$$

*the Hankel determinants satisfy*

$$\det\left( f_{i+j-m} \right)_{i,j=0}^{n} = 0 \quad \text{for } n < m \tag{2.33}$$

$$\det\left( f_{i+j-m} \right)_{i,j=0}^{m} = (-1)^{\binom{m+1}{2}} \tag{2.34}$$

*and*

$$\det\left( f_{i+j-m} \right)_{i,j=0}^{n+m} = (-1)^{\binom{m+1}{2}} a^n \det\left( g_{i+j+m-p+2} \right)_{i,j=0}^{n-1} \tag{2.35}$$

*for* $n > 0$ *if we define* $g_n = 0$ *for* $-p < n < 0$.

## 3. The main results

**Theorem 3.1 (V. I. Buslaev [4])**

*Let* $(b_n)_{n\geq -1}$ *be a non-decreasing sequence of integers such that* $b_{n+2} - b_n \geq 1$ *with initial values* $b_{-1} = -1$ *and* $b_0 = 0$ *and let*

$$f(x) = \sum_{n\geq 0} f_n x^n = \cfrac{1}{1 - \cfrac{a_0 x^{b_1 - b_{-1}}}{1 - \cfrac{a_1 x^{b_2 - b_0}}{1 - \cfrac{a_2 x^{b_3 - b_1}}{1 - \cfrac{a_3 x^{b_4 - b_2}}{1 - \ddots}}}}}. \tag{3.1}$$

*Then the Hankel determinants*

$$d(n) = \det\left( f_{i+j} \right)_{i,j=0}^{n} \tag{3.2}$$

*satisfy*

$$d(b_k) = (-1)^{\binom{b_1 - b_0}{2} + \binom{b_2 - b_1}{2} + \cdots + \binom{b_k - b_{k-1}}{2}} a_0^{b_k - b_0} a_1^{b_k - b_1} a_2^{b_k - b_2} a_3^{b_k - b_3} \cdots a_{k-1}^{b_k - b_{k-1}} \tag{3.3}$$

*and vanish for all other values of* $n > 0$.

**Proof**

Let

$$f^{(k)}(x) = \frac{1}{1 - a_k x^{b_{k+1}-b_{k-1}} f^{(k+1)}(x)}. \tag{3.4}$$

Then $f(x) = f^{(0)}(x).$

We want to prove that

$$\det\left(f_{i+j}^{(0)}\right)_{i,j=0}^{n+b_k} = (-1)^{\sum_{j=0}^{k-1}\binom{b_{j+1}-b_j}{2}} \prod_{j=0}^{k} a_j^{n+b_k-b_j} \det\left(f_{i+j+1+b_k-b_{k+1}}^{(k+1)}\right)_{i,j=0}^{n-1} \tag{3.5}$$

if we set $\det\left(f_{i+j+1+b_k-b_{k+1}}^{(k+2)}\right)_{i,j=0}^{n-1} = 1$ for $n = 0.$

This is true for $k = 0.$ In this case it reduces to

$$\det\left(f_{i+j}^{(0)}\right)_{i,j=0}^{n} = a_0^n \det\left(f_{i+j+1-b_1}^{(1)}\right)_{i,j=0}^{n-1}.$$

Now suppose that (3.5) is true for $k.$ This means that

$$\det\left(f_{i+j}^{(0)}\right)_{i,j=0}^{n+b_k} = (-1)^{\sum_{j=0}^{k-1}\binom{b_{j+1}-b_j}{2}} \prod_{j=0}^{k} a_j^{n+b_k-b_j} \det\left(f_{i+j+1+b_k-b_{k+1}}^{(k+1)}\right)_{i,j=0}^{n-1}.$$

Therefore

$$\det\left(f_{i+j}^{(0)}\right)_{i,j=0}^{n+b_{k+1}} = (-1)^{\sum_{j=0}^{k-1}\binom{b_{j+1}-b_j}{2}} \prod_{j=0}^{k} a_j^{n+b_{k+1}-b_j} \det\left(f_{i+j+1+b_k-b_{k+1}}^{(k+1)}\right)_{i,j=0}^{n-1+b_{k+1}-b_k}.$$

By Lemma 2.1 we have

$$\det\left(f_{i+j+1+b_k-b_{k+1}}^{(k+1)}\right)_{i,j=0}^{n-1+b_{k+1}-b_k} = (-1)^{\binom{b_{k+1}-b_k}{2}} a_{k+1}^n \det\left(f_{i+j+1+b_{k+1}-b_{k+2}}^{(k+2)}\right)_{i,j=0}^{n-1}.$$

Therefore (3.5) is proved.

Theorem 3.1 follows since $\det\left(f_{i+j+1+b_{k+1}-b_{k+2}}^{(k+2)}\right)_{i,j=0}^{n-1} = 0$ for $0 < n < b_{k+2} - b_{k+1} - 1$.

**Remark**

This theorem has also been proved by V.I. Buslaev [4] with another method.

It would also be interesting to have formulae for the Hankel determinants of

$$f(x)=\sum_{n\geq 0} f_n x^n = \cfrac{1}{1-\cfrac{a_0x^{m_0}}{1-\cfrac{a_1x^{m_1}}{1-\cfrac{a_2x^{m_2}}{1-\cfrac{a_3x^{m_3}}{1-\ddots}}}}}$$

for an arbitrary sequence of positive integers $m_n$ instead of $b_n - b_{n-2}$. This has been attempted in [3], where it is claimed that all such Hankel determinants are products of $a_n$'s. But there are simple counter examples.

Let for example $(m_n)=(1,2,1,1,1,\cdots)$. Then the Hankel determinants are

$$1,\ 0,\ -(a_0a_1)^2,\ -(a_0a_1)^3(a_2a_3)(a_3+a_4),\ -(a_0a_1)^4(a_2a_3)^2(a_4a_5)(a_3a_5+a_3a_6+a_4a_6),\cdots.$$

The general terms will be given in (3.19). The corresponding sequence $(b_n)_{n\geq 0}=(0,0,2,1,3,2,4,3,5,4,\cdots)$ does not satisfy our assumptions.

If each $a_k=\pm 1$ in (3.1) then each non-zero Hankel determinant is also $\pm 1$.

Note that $m_n = b_{n+1}-b_{n-1}$ implies $b_{2n}=\sum_{j=0}^{n-1} m_{2n-1-2j}$ and $b_{2n+1}=\sum_{j=0}^{n} m_{2n-2j}-1$.

The set $\{b_n\}_{n\geq 0}$ is the set of all $n\in\mathbb{N}$ such that $d(n)\neq 0$. If $b_i = b_{i+1}$ then $d(b_k)$ contains the term $(a_ia_{i+1})^{b_k-b_i}$, i.e. the powers of $a_j$ in $d(n)$ have a similar "pattern" as the sequence $(b_n)_{n\geq 0}$.

Let for example $(b_n)_{n\geq 0}=(0,0,1,2,5,6,6,8,9,10,\cdots)$ i.e. $(p_n)_{n\geq 0}=(1,1,2,4,4,1,2,3,2,2,2,\cdots)$.

Then the sequence of Hankel determinants is

$$1,\ (a_0a_1),\ (a_0a_1)^2a_2,\ 0,\ 0,\ -(a_0a_1)^5a_2^4a_3^3,\ -(a_0a_1)^6a_2^5a_3^4a_4,\ 0,\ (a_0a_1)^8a_2^7a_3^6a_4^3(a_5a_6)^2,$$
$$(a_0a_1)^9a_2^8a_3^7a_4^4(a_5a_6)^3a_7,\ \cdots.$$

**Theorem 3.2**

*Let $(b_k)_{k\geq -1}$ be a non-decreasing sequence of integers such that $b_{k+2}-b_k\geq 1$ with initial values $b_{-1}=-1$ and $b_0=0$ and let $(B_n)_{n\geq 0}$ be the elements of the set $\{b_k\}_{k\geq -1}$ in increasing order.*

*If $B_n = b_{k-1}$ has multiplicity 1 or if $B_n=b_{k-2}=b_{k-1}$ has multiplicity 2 then $P_n(x):=\dfrac{p_{B_n+1}(x)}{d(B_n)}=r_k(x)$.*

*For $B_n+1<m<B_{n+1}-1$ we get $p_m(x)=0$ and for $m=B_{n+1}-1$ we get*

$$p_{B_{n+1}}(x)=\pm(a_0a_1\cdots a_{k-1})^{B_{n+1}-B_n-1}p_{B_n+1}(x).$$

**Proof**

Consider

$$\Lambda\left(\frac{1}{\det(f_{i+j})_{i,j=0}^{b_{k-1}}}\det\begin{pmatrix} f_0 & f_1 & \cdots & f_{b_{k-1}} & 1 \\ f_1 & f_2 & \cdots & f_{b_{k-1}+1} & x \\ f_2 & f_3 & \cdots & f_{b_{k-1}+2} & x^2 \\ \vdots & & & & \vdots \\ f_{b_{k-1}+1} & f_{b_{k-1}+2} & \cdots & f_{2b_{k-1}+1} & x^{b_{k-1}+1} \end{pmatrix} x^\ell\right)$$

$$=\frac{1}{\det(f_{i+j})_{i,j=0}^{b_{k-1}}}\det\begin{pmatrix} f_0 & f_1 & \cdots & f_{b_{k-1}} & f_\ell \\ f_1 & f_2 & \cdots & f_{b_{k-1}+1} & f_{\ell+1} \\ f_2 & f_3 & \cdots & f_{b_{k-1}+2} & f_{\ell+2} \\ \vdots & & & & \vdots \\ f_{b_{k-1}+1} & f_{b_{k-1}+2} & \cdots & f_{2b_{k-1}+1} & f_{\ell+b_{k-1}+1} \end{pmatrix}$$

By (2.12) we have $\sum_j v_{k,j} f_{b_{k-1}-j+\ell+1}=0$ for $\ell<b_k$. This means that the last column can be reduced to the $0-$column by elementary column operations. Therefore the determinant vanishes for these $\ell$.
For $\ell=b_k$ we get by (2.13)

$\sum_j v_{k,j} f_{b_{k-1}-j+b_k+1}=a_0\cdots a_{k-1}$ and therefore by elementary column operations

$$\Lambda\left(\frac{p_{B_n+1}(x)}{d(B_n)}x^{b_k}\right)=\frac{1}{\det(f_{i+j})_{i,j=0}^{b_{k-1}}}\det\begin{pmatrix} f_0 & f_1 & \cdots & f_{b_{k-1}} & 0 \\ f_1 & f_2 & \cdots & f_{b_{k-1}+1} & 0 \\ f_2 & f_3 & \cdots & f_{b_{k-1}+2} & 0 \\ \vdots & & & & \vdots \\ f_{b_{k-1}+1} & f_{b_{k-1}+2} & \cdots & f_{2b_{k-1}+1} & a_0\cdots a_{k-1} \end{pmatrix}=a_0\cdots a_{k-1}.$$

Therefore the polynomials $q(x):=\frac{p_{B_n+1}(x)}{d(B_n)}$ which has degree $b_{k-1}+1$ and $r_k(x)$ which is of the same degree satisfy $\Lambda\left(\left(q(x)-r_k(x)\right)x^n\right)=0$ for $0\le n\le b_k$. This implies $q(x)=r_k(x)$ since the rank of the matrix $\begin{pmatrix} f_0 & f_1 & \cdots & f_{b_{k-1}} & f_{b_k} \\ f_1 & f_2 & \cdots & f_{b_{k-1}+1} & f_{b_k+1} \\ f_2 & f_3 & \cdots & f_{b_{k-1}+2} & f_{b_k+2} \\ \vdots & & & & \vdots \\ f_{b_{k-1}+1} & f_{b_{k-1}+2} & \cdots & f_{2b_{k-1}+1} & f_{b_k+b_{k-1}+1} \end{pmatrix}$ is $b_{k-1}+1$.

To prove the other assertion observe that (2.12) implies that at least one column in (2.14) can be reduced to the zero column and therefore $p_m(x)=0$ for $B_n+1<m<B_{n+1}-1.$

For $m=B_{n+1}=b_k$ column reduction reduces column $b_{k-1}+1$
which is

$$\begin{pmatrix} f_{b_{k-1}+1} \\ f_{b_{k-1}+2} \\ \vdots \\ f_{b_k-1} \\ f_{b_k} \end{pmatrix} \text{ to } \begin{pmatrix} 0 \\ \vdots \\ 0 \\ a_0\cdots a_{k-1} \end{pmatrix}.$$

Now expand the determinant with respect to this column and iterate $b_k-b_{k-1}-1$ times.

As an example consider $\left(b_k\right)_{k\ge -1}=\left(-1,0,0,3,3,7,8,9,\cdots\right)$ and therefore $\left(B_n+1\right)_{n\ge 0}=\left(0,1,4,8,9,10,\cdots\right).$
For $B_{-1}=b_{-1}=-1$ we set $p_0(x)=r_0(x)=1$.
For $B_0=0=b_0=b_1$ we get $\dfrac{p_1(x)}{d(0)}=r_2(x),$
for $B_1=3=b_2=b_3$ we get $\dfrac{p_4(x)}{d(3)}=r_4(x),$
for $B_2=7=b_4$ we get $\dfrac{p_8(x)}{d(7)}=r_5(x),$ etc.

The remaining polynomials $p_n(x)$ are $p_2(x)=0,$
$p_3(x)=-(a_0a_1)^2r_2,$ because $p_3(x)=p_{B_2}(x)$ and $B_1=b_0=b_1$,
$p_5(x)=p_6(x)=0,$
$p_7(x)=p_{B_3}(x)$ and since $B_2=b_2=b_3$ we have
$p_7(x)=p_{B_3}(x)=\left(a_0a_1a_2a_3\right)^3p_4(x)=\left(a_0a_1a_2a_3\right)^3d(3)r_4(x)$, etc.
Let us show in detail how to derive $p_3(x).$

$$p_3(x)=\det\begin{pmatrix} 1 & a_0 & a_0^2 & 1 \\ a_0 & a_0^2 & a_0^3 & x \\ a_0^2 & a_0^3 & a_0^4+a_0a_1 & x^2 \\ a_0^3 & a_0^4+a_0a_1 & a_0^5+2a_0^2a_1 & x^3 \end{pmatrix}=\det\begin{pmatrix} 1 & 0 & 0 & 1 \\ a_0 & 0 & 0 & x \\ a_0^2 & 0 & a_0a_1 & x^2 \\ a_0^3 & a_0a_1 & a_0a_1 & x^3 \end{pmatrix}=-(a_0a_1)^2\det\begin{pmatrix} 1 & 1 \\ a_0 & x \end{pmatrix}.$$

First we give two well-known examples from this point of view.

**Example 3.1**

The special case $b_n = n$ and therefore $m_n = 2$ gives again (2.21) and therefore
$d(n) = a_0^n a_1^{n-1} a_2^{n-2} \cdots a_{n-1}.$
The polynomials $r_k(x)$ have degree $\deg r_k = k$ and satisfy $r_k(x) = x r_{k-1}(x) - a_{k-2} r_{k-2}(x)$ with $r_{-1}(x) = 0$ and $r_0(x) = 1.$ They are orthogonal with respect to $\Lambda$ and satisfy $r_k(x) = \frac{p_k(x)}{d(k-1)} = P_k(x).$
By (2.14) we get $d_1(k) = \det\left( f_{i+j+1} \right)_{i,j=0}^{n} = p_{k+1}(0)$ and therefore $\frac{d_1(k)}{d(k)} = r_{k+1}(0).$
It is clear that $r_{2k+1}(0) = 0$ and therefore $d_1(2n) = 0$.
Furthermore we get $r_{2k}(0) = -a_{2k-2} r_{2k-2}(0) = \cdots = (-1)^{k-1} a_0 a_2 \cdots a_{2(k-1)}.$
This gives for $n > 0$

$$d_1(2n-1) = (-1)^n a_0^{2n} \left(a_1 a_2\right)^{2n-2} \left(a_3 a_4\right)^{2n-4} \cdots \left(a_{2n-3} a_{2n-2}\right)^2. \tag{3.6}$$

The most important special case and prototype for many papers on Hankel determinants is

**Example 3.1.1**

Let $b_n = n$ and all $a_n = 1.$ Then by (2.24) $f(x) = \sum_{n \geq 0} C_n x^{2n}.$ In this case $d(n) = 1$ for all $n.$ By this property the Catalan numbers are uniquely determined.

If we define the Fibonacci polynomials $Fib_n(x) = \sum_{k=0}^{\left\lfloor \frac{n-1}{2} \right\rfloor} \binom{n-1-k}{k} (-1)^k x^{n-1-2k}$ by
$Fib_n(x) = x Fib_{n-1}(x) - Fib_{n-2}(x)$ with initial values $Fib_0(x) = 0$ and $Fib_1(x) = 1$, then

$$r_k(x) = \det\left( f_{i+j} x - f_{i+j+1} \right)_{i,j=0}^{k-1} = Fib_{k+1}(x) \tag{3.7}$$

if $f_{2n} = C_n$ and $f_{2n+1} = 0.$

In this case (2.12) and (2.13) are equivalent with
$\Lambda\left(Fib_{2k}(x) x^{2n}\right) = \sum_{j=0}^{k} (-1)^j \binom{2k-j}{j} C_{k+n-j} = 0$ for $n < k$ and $= 1$ for $n = k,$
and to

$\Lambda\left(Fib_{2k+1}(x) x^{2n+1}\right) = \sum_{j=0}^{k} (-1)^j \binom{2k+1-j}{j} C_{k+n+1-j} = 0$ for $n < k$ and $= 1$ for $n = k.$

**Example 3.2**

For $(b_n)_{n\geq 0} = (0,0,1,1,2,2,3,3,\cdots)$ we get

$$f(x) = \sum_{n\geq 0} f_n x^n = \cfrac{1}{1-\cfrac{a_0 x}{1-\cfrac{a_1 x}{1-\cfrac{a_2 x}{1-\ddots}}}}. \tag{3.8}$$

and therefore as in (2.31)

$$d(n) = (a_0 a_1)^n (a_2 a_3)^{n-1} \cdots (a_{2n-2} a_{2n-1}). \tag{3.9}$$

In this case

$$d_1(n) = a_0^{n+1} (a_1 a_2)^n (a_3 a_4)^{n-1} \cdots (a_{2n-1} a_{2n}). \tag{3.10}$$

This follows immediately from $\det\left(f_{i+j+1}\right)_{i,j=0}^{n} = a_0^{n+1} \det\left(f_{i+j}^{(1)}\right)_{i,j=0}^{n}$ and (3.9).

Later we shall need $d_2(n) = \det\left(f_{i+j+2}\right)_{i,j=0}^{n}$. To compute this determinant we can use the condensation formula (cf. [8], (2.16))

$$d_2(n)d(n) = d_2(n-1)d(n+1) + d_1(n)^2, \tag{3.11}$$

which by (3.9) and (3.10) reduces to

$$\frac{d_2(n)}{d_1(n)} = a_{2n+1} \frac{d_2(n-1)}{d_1(n-1)} + a_0 a_2 \cdots a_{2n} = \sum_{i_1 < i_2 - 1 < i_3 - 2 < \cdots < i_n - n + 1 \leq n} a_{i_1} a_{i_2} \cdots a_{i_n}.$$

Therefore we get

$$d_2(n) = a_0^{n+1} (a_1 a_2)^n (a_3 a_4)^{n-1} \cdots (a_{2n-1} a_{2n}) \sum_{i_1 < i_2 - 1 < i_3 - 2 < \cdots < i_n - n + 1 \leq n} a_{i_1} a_{i_2} \cdots a_{i_n}. \tag{3.12}$$

If (3.8) holds we have $b_{2k} = b_{2k+1} = k$ and $B_n = n-1.$ Thus $B_n = n-1 = b_{2n-2} = b_{2n-1}$ and
$r_{2n}(x) = \dfrac{p_n(x)}{d(n-1)} = P_n(x).$

In this case there is also a simple recursion for $P_n(x)$.

We have $P_1(x) = x - a_0$ and $P_n(x) = (x - a_{2n-2} - a_{2n-3}) P_{n-1}(x) - a_{2n-4} a_{2n-3} P_{n-2}(x)$ for $n \geq 2.$

This follows from

$$r_{2n}(x) = r_{2n-1} - a_{2n-2} r_{2n-2} = x r_{2n-2}(x) - a_{2n-2} r_{2n-2} - a_{2n-3} r_{2n-3}$$
$$= x r_{2n-2}(x) - a_{2n-2} r_{2n-2}(x) - a_{2n-3}\left(r_{2n-2}(x) + a_{2n-4} r_{2n-4}(x)\right).$$

**Corollary 3.1**

*A series* $f(x)=\sum_{n\geq 0} f_n x^n$ *has a representation as a continued fraction of the form* (3.8) if and only if $d(n)\neq 0$ and $d_1(n)\neq 0$ *for all* $n\in\mathbb{N}$.

**Proof**

If no determinant vanishes then the numbers $a_n$ can be uniquely computed from (3.9) and (3.10). This is clear from

$\frac{d_1(n)}{d(n)}=a_0a_2\cdots a_{2n}$ and $\frac{d(n)}{d_1(n-1)}=a_1a_3\cdots a_{2n-1}$.

**Example 3.2.1**

The simplest example is $f(x)=\sum_{n\geq 0} C_n x^n$ where $d(n)=d_1(n)=1$.

Since in this case $d_2(n)=n+2$ by (3.11) we see that $f(x)=\sum_{n\geq 0} C_{n+1}x^n$ has a representation of the form (3.8) with $a_{2n}=\frac{n+2}{n+1}$ and $a_{2n+1}=\frac{n+1}{n+2}$.

**Example 3.2.2**

G. Eisenstein has shown (cf. Perron [9], §59, (32)) that $f(x)=\sum_{n\geq 0} q^{\binom{n+1}{2}}x^n$ has an expansion of the form (3.8) with $a_{2n}=q^{2n+1}$ and $a_{2n+1}=\left(q^{n+1}-1\right)q^{n+1}$. (A simple proof can also be found in [6], (3.3) and (3.6)). For $q\to 1$ we get $a_{2n}\to 1$ and $a_{2n+1}\to 0$. So in some sense the continued fractions converge to $\frac{1}{1-x}=\sum_{n\geq 0} x^n$.

The Hankel determinants are $d(n)=q^{\frac{n(n+1)^2}{2}}\prod_{j=1}^{n}\left(q^j-1\right)^{n+1-j}$.

The polynomials $r_k(x)$ are given by $r_{2k}(x)=\sum_{j=0}^{k}(-1)^j q^{kj}\begin{bmatrix}k\\j\end{bmatrix}x^{k-j}$ and

$r_{2k+1}(x)=\sum_{j=0}^{k}(-1)^j\begin{bmatrix}k\\j\end{bmatrix}q^{j(k+1)}x^{k+1-j}$, where $\begin{bmatrix}n\\k\end{bmatrix}=\prod_{j=0}^{k-1}\frac{1-q^{n-j}}{1-q^{k-j}}$ is a $q-$binomial coefficient.

These assertions are by (2.12) and (2.13) equivalent with

$$\Lambda\left(r_{2k}(x)x^m\right)=\sum_{j=0}^{k}(-1)^j q^{kj}\begin{bmatrix}k\\j\end{bmatrix}q^{\binom{k-j+m+1}{2}}=0 \text{ for } m<k \text{ and}$$

$$\Lambda\left(r_{2k}(x)x^k\right)=\sum_{j=0}^{k}(-1)^j q^{kj}\begin{bmatrix}k\\j\end{bmatrix}q^{\binom{2k-j+1}{2}}=q^{k^2+\binom{k+1}{2}}(q-1)(q^2-1)\cdots(q^k-1)$$

and

$$\Lambda\left(r_{2k+1}(x)x^m\right)=\sum_{j=0}^{k}(-1)^j q^{(k+1)j}\begin{bmatrix}k\\j\end{bmatrix}q^{\binom{k-j+m+2}{2}}=0 \text{ for } k<m \text{ and}$$

$$\Lambda\left(r_{2k+1}(x)x^k\right)=\sum_{j=0}^{k}(-1)^j q^{(k+1)j}\begin{bmatrix}k\\j\end{bmatrix}q^{\binom{2k-j+2}{2}}=q^{(k+1)^2+\binom{k+1}{2}}(q-1)(q^2-1)\cdots(q^k-1).$$

It is easy to show these identities directly by using the well-known formula

$$\sum_{j=0}^{k}\begin{bmatrix}k\\j\end{bmatrix}q^{\binom{j}{2}}z^j=(1+z)(1+qz)\cdots(1+q^{k-1}z).$$

We get

$$\sum_{j=0}^{k}(-1)^j q^{kj}\begin{bmatrix}k\\j\end{bmatrix}q^{\binom{k-j+m+1}{2}}=q^{\binom{k+1}{2}+\binom{m+1}{2}+km}\sum_{j=0}^{k}(-1)^j\begin{bmatrix}k\\j\end{bmatrix}q^{\binom{j}{2}-jm}=q^{\binom{k+1}{2}+\binom{m+1}{2}+km}\left(1-q^{-m}\right)\left(1-q^{1-m}\right)\cdots\left(1-q^{k-1-m}\right).$$

For $m<k$ the right-hand side vanishes and for $m=k$ get $q^{\binom{k+1}{2}+k^2}(q-1)\cdots(q^k-1)$

For $r_{2k+1}(x)$ we get

$$\sum_{j=0}^{k}(-1)^j q^{(k+1)j}\begin{bmatrix}k\\j\end{bmatrix}q^{\binom{k-j+m+2}{2}}=q^{\binom{k+2}{2}+\binom{m+2}{2}+km-1}\sum_{j=0}^{k}(-1)^j\begin{bmatrix}k\\j\end{bmatrix}q^{\binom{j}{2}-jm}$$

$$=q^{\binom{k+2}{2}+\binom{m+2}{2}+km-1}\left(1-q^{-m}\right)\left(1-q^{1-m}\right)\cdots\left(1-q^{k-1-m}\right).$$

For $m=k$ this gives $q^{\binom{k+1}{2}+(k+1)^2}(q-1)\cdots(q^k-1)$.

As a special case we get that

$$\frac{\det\left(q^{\binom{i+j+1}{2}}x-q^{\binom{i+j+2}{2}}\right)_{i,j=0}^{k-1}}{d(k-1)}=r_{2k}(x)=\sum_{j=0}^{k}(-1)^j q^{kj}\begin{bmatrix}k\\j\end{bmatrix}x^{k-j}.$$

For $q\to 1$ these polynomials converge to $r_{2k}(x)=(x-1)^k$ and $r_{2k+1}(x)=x(x-1)^k$, which are no longer in an analogous relation to the finite continued fraction $\frac{1}{1-x}$.

**Example 3.2.3**

For the sequence of Motzkin numbers $(M_n)_{n\geq 0} = (1,1,2,4,9,21,51,\cdots)$ whose generating function $f(x) = \sum_{n\geq 0} M_n x^n = \frac{1-x-\sqrt{1-2x-3x^2}}{2x^2}$ satisfies

$$f(x) = 1 + xf(x) + x^2 f(x)^2 \tag{3.13}$$

the Hankel determinants are $d(n) = 1$ and $(d_1(n))_{n\geq 0} = (1,0,-1,-1,0,1,1,0,-1,-1,0,1,\cdots)$ which is periodic with period $6.$
Since some determinants vanish there is no representation of the form (3.8).

We have instead

$$f(x) = \cfrac{1}{1-\cfrac{x}{1-\cfrac{x}{1-x^2 f(x)}}}.$$

For (3.13) implies $\frac{1}{f(x)} - 1 = -x(1+xf(x))$ and $\frac{1}{1+xf(x)} = 1 - \frac{x}{1-x^2 f(x)}.$
Thus in this case $(b_n)_{n\geq -1} = (-1,0,0,1,2,2,3,4,4,5,6,6,7,\cdots)$ and $(m_n)_{n\geq 0} = (1,1,2,1,1,2,\cdots).$

**Example 3.2.4**

If we consider the numbers $M_n(u) = \sum_{k=0}^{\left\lfloor \frac{n}{2} \right\rfloor} \binom{2k}{k} \frac{1}{k+1} \binom{n}{2k} u^{n-2k}$ whose generating function $f(x,u)$ satisfies $f(x,u) = 1 + uxf(x,u) + x^2 f(x,u)^2$ then the determinants are $d(n) = 1$ and $d_1(n) = Fib_{n+2}(u)$ (cf. e.g. [5],Theorem A*).
The sequence $(Fib_n(1)) = (0,1,1,0,-1,-1,\cdots)$ is periodic with period $6.$
For $u \to 1$ the numbers $M_n(u)$ converge to the Motzkin numbers $M_n.$
Therefore for $u \neq 1$ in some neighborhood of $1$ all $f(x,u)$ have a representation of the form (3.8).

It is given by $a_{2n} = \frac{Fib_{n+2}(u)}{Fib_{n+1}(u)}$ and $a_{2n-1} = \frac{Fib_n(u)}{Fib_{n+1}(u)}.$ For $u \to 1$ we have $a_0 \to 1,$ $a_1 \to 1,$ and $a_2 \to 0.$ This explains why we have for the Motzkin numbers $m_2 > 1.$

**Remark**
If the sequence $(b_n)_{n\geq 0}$ contains all non-negative integers then $d(n) \neq 0$ for all $n.$ In this case the above results are special cases of the well-known connection with orthogonal polynomials and continued J-fractions (cf. e.g. [8], (2.30) or [5]). Observe also that in this case all $m_n \in \{1,2\}$ and each sequence of $1'\text{s}$ must have even length.

**Example 3.2.5**

As in the example with the Motzkin numbers let the power sequence be $(m_n)_{n\geq 0} = (1,1,2,1,1,2,\cdots)$, i.e. $(b_n)_{n\geq -1} = (-1,0,0,1,2,2,3,4,4,5,6,6,7,\cdots)$ with $b_{3k} = b_{3k+1} = 2k$ and $b_{3k+2} = 2k+1.$
Since $B_{2n} = 2n-1 = b_{3n-1}$ and $B_{2n+1} = 2n = b_{3n} = b_{3n+1}$ we get

$$P_{2n}(x) = \frac{p_{B_{2n}+1}(x)}{d(B_{2n})} = r_{3n}(x)$$

and

$$P_{2n+1}(x) = \frac{p_{B_{2n+1}+1}(x)}{d(B_{2n+1})} = r_{3n+2}(x).$$

The recurrences are

$$P_{2n+1}(x) = r_{3n+1}(x) - a_{3n}r_{3n}(x) = \left(xr_{3n}(x) - a_{3n-1}r_{3n-1}(x)\right) - a_{3n}r_{3n}(x) = (x-a_{3n})P_{2n}(x) - a_{3n-1}P_{2n-1}(x)$$

and since $r_{3n-1}(x) = r_{3n-2}(x) - a_{3n-3}r_{3n-3}(x)$
we get

$$P_{2n}(x) = r_{3n}(x) = xr_{3n-1}(x) - a_{3n-2}r_{3n-2}(x) = xr_{3n-1}(x) - a_{3n-2}\left(r_{3n-1}(x) + a_{3n-3}r_{3n-3}(x)\right)$$
$$= (x-a_{3n-2})P_{2n-1}(x) - a_{3n-3}a_{3n-2}P_{2n-2}(x).$$

By Theorem 3.1 the Hankel transform is
$1,\ (a_0a_1),\ (a_0a_1)^2 a_2,\ (a_0a_1)^3 a_2^2 (a_3a_4),\ (a_0a_1)^4 a_2^3 (a_3a_4)^2 a_5, \cdots$
or
$d(2n) = d(b_{3n}) = d(b_{3n+1}) = (a_0a_1)^{2n} a_2^{2n-1} (a_3a_4)^{2n-2} \cdots a_{3n-1}$
and
$d(2n+1) = d(b_{3n+2}) = (a_0a_1)^{2n+1} a_2^{2n} (a_3a_4)^{2n-1} \cdots (a_{3n}a_{3n+1}).$
Since $P_2(x) = (x-a_1)(x-a_0) - a_0a_1 = x(x-a_0-a_1)$ we see that $P_2(0)=0$ and therefore $d_1(1)=0.$
The transformation into a J-fraction gives

$$\cfrac{1}{1-\cfrac{a_0x}{1-\cfrac{a_1x}{1-\cfrac{a_2x^2}{1-\cfrac{a_3x}{1-\cfrac{a_4x}{1-\cfrac{a_5x^2}{1-\ddots}}}}}}} = \cfrac{1}{1-a_0x-\cfrac{a_0a_1x^2}{1-a_1x-\cfrac{a_2x^2}{1-a_3x-\cfrac{a_3a_4x^2}{1-a_4x-\cfrac{a_5x^2}{1-a_6x-\ddots}}}}}$$

**Example 3.3**

Let $(b_k)_{k\geq -1} = (-1,0,0,3,5,5,8,10,10,13,\cdots)$ with $b_{3k} = b_{3k+1} = 5k$ and $b_{3k+2} = 5k+3.$

Here we also have $P_{2n}(x) = r_{3n}(x)$ and $P_{2n+1}(x) = r_{3n+2}(x).$
We get
$$P_{2n+1}(x) = r_{3n+1}(x) - a_{3n}r_{3n}(x) = \left(x^2 r_{3n}(x) - a_{3n-1}r_{3n-1}(x)\right) - a_{3n}r_{3n}(x) = (x^2 - a_{3n})P_{2n}(x) - a_{3n-1}P_{2n-1}(x)$$

and from $r_{3n-1}(x) = r_{3n-2}(x) - a_{3n-3}r_{3n-3}(x)$
we get

$$P_{2n}(x) = r_{3n}(x) = x^3 r_{3n-1}(x) - a_{3n-2}r_{3n-2}(x) = x^3 r_{3n-1}(x) - a_{3n-2}\left(r_{3n-1}(x) + a_{3n-3}r_{3n-3}(x)\right)$$
$$= (x^3 - a_{3n-2})P_{2n-1}(x) - a_{3n-3}a_{3n-2}P_{2n-2}(x).$$

**Example 3.4**

Let $m_n = m \geq 1$ for all $n.$ This is equivalent with $b_{2k} = mk$ and $b_{2k-1} = mk - 1.$
The corresponding polynomials $r_k(x) = P_k^{(m)}(x)$ are given by $r_0(x) = 1,$ $r_1(x) = x,$
$$r_{2k}(x) = x^{m-1}r_{2k-1}(x) - a_{2k-2}r_{2k-2}(x)$$
and
$$r_{2k+1}(x) = xr_{2k}(x) - a_{2k-1}r_{2k-1}(x).$$

Then $\deg P_{2k}^{(m)}(x) = mk$ and $\deg P_{2k+1}^{(m)}(x) = mk+1.$

Let
$$f_2(x) = \cfrac{1}{1 - \cfrac{a_0x^2}{1 - \cfrac{a_1x^2}{1 - \cfrac{a_2x^2}{1 - \ddots}}}} = \sum_{n\geq 0} c_n x^{2n}.$$
Then we have
$$f_m(x) = \cfrac{1}{1 - \cfrac{a_0x^m}{1 - \cfrac{a_1x^m}{1 - \cfrac{a_2x^m}{1 - \ddots}}}} = \sum_{n\geq 0} c_n x^{mn}.$$

**Proposition 3.1**
*Let* $P_{2k}^{(2)}(x) = \sum_{j=0}^{k} u_{k,j}x^{2j}$ *and* $P_{2k-1}^{(2)}(x) = x\sum_{j=0}^{k-1} v_{k,j}x^{2j}.$
*Then*
$$P_{2k}^{(m)}(x) = \sum_{j=0}^{k} u_{k,j}x^{mj} \text{ and } P_{2k-1}^{(m)}(x) = x\sum_{j=0}^{k-1} v_{k,j}x^{mj}.$$

**Proof**

Define linear functionals $\Lambda_m$ by $\Lambda_m\left(x^{mn}\right)=c_n$ and $\Lambda_m\left(x^{mn+j}\right)=0$ for $0<j<m.$

The polynomials $P_k^{(m)}(x)$ are uniquely determined by

$$\Lambda_m\left(P_{2k}^{(m)}(x)x^n\right)=0 \text{ for } n<mk \text{ and } \Lambda_m\left(P_{2k}^{(m)}(x)x^{mk}\right)=a_0\cdots a_{2k-1}$$

and

$$\Lambda_m\left(P_{2k-1}^{(m)}(x)x^n\right)=0 \text{ for } n<mk-1 \text{ and } \Lambda_m\left(P_{2k-1}^{(m)}(x)x^{mk-1}\right)=a_0\cdots a_{2k-2}.$$

Now we have

$$\Lambda_m\left(x^{mi}\sum_{j=0}^{k}u_{k,j}x^{mj}\right)=\Lambda_2\left(x^{2i}\sum_{j=0}^{k}u_{k,j}x^{2j}\right)=0 \text{ for } i<k \text{ and } \Lambda_m\left(x^{n}\sum_{j=0}^{k}u_{k,j}x^{mj}\right)=0 \text{ if } n<mk \text{ and }$$

$n\not\equiv 0 \bmod m$

and

$$\Lambda_m\left(x^{mk}\sum_{j=0}^{k}u_{k,j}x^{mj}\right)=\Lambda_2\left(x^{2k}\sum_{j=0}^{k}u_{k,j}x^{2j}\right)=a_0\cdots a_{k-1}.$$

Therefore we must have $P_{2k}^{(m)}(x)=\sum_{j=0}^{k}u_{k,j}x^{mj}$ .

The same argument applies in the second case.

For the Hankel determinants we get

$$d(nm)=d\left(b_{2n}\right)=(-1)^{\binom{m-1}{2}n}a_0^{nm}a_1^{1+(n-1)m}a_2^{(n-1)m}a_3^{1+(n-2)m}\cdots a_{2n-2}^{m}a_{2n-1},$$

$$d(nm-1)=d\left(b_{2n-1}\right)=(-1)^{\binom{m-1}{2}n}a_0^{nm-1}a_1^{(n-1)m}a_2^{(n-1)m-1}\cdots a_{2n-3}^{m}a_{2n-2}^{m-1}.$$

If all $a_i=a$ we get $f_m(x)=\sum_n f_n^{(m)}x^n=\sum_n a^nC_nx^{mn},$ where $C_n=\binom{2n}{n}\frac{1}{n+1}$ are the Catalan numbers.

The corresponding Hankel determinants which have been computed in [5] with another method are

$$\det\left(f_{i+j}^{(m)}\right)_{i,j=0}^{mn-1}=(-1)^{\binom{m-1}{2}n}a^{n(mn-1)},$$
$$\det\left(f_{i+j}^{(m)}\right)_{i,j=0}^{mn}=(-1)^{\binom{m-1}{2}n}a^{n(mn+1)}, \tag{3.14}$$

and

$\det\left(f_{i+j}^{(m)}\right)_{i,j=0}^{n}=0$ for all other $n.$

Under the same assumptions we have

$$\det\left(f_{i+j+1}^{(m)}\right)_{i,j=0}^{mn-1} = (-1)^{\binom{m}{2}n} a^{mn^2}. \tag{3.15}$$

Since

$$\det\left(f_{i+j+1}^{(m)}\right)_{i,j=0}^{n-1} = a^n \det\left(f_{i+j-m+1}^{(m)}\right)_{i,j=0}^{n-1}$$

we get $\det\left(f_{i+j+1}^{(m)}\right)_{i,j=0}^{n-1} = 0$ for $0 < n < m-1$ and

$$\det\left(f_{i+j+1}^{(m)}\right)_{i,j=0}^{n+m-1} = a^{n+m} \det\left(f_{i+j-m+1}^{(m)}\right)_{i,j=0}^{n+m-1} = (-1)^{\binom{m}{2}} a^{2n+m} \det\left(f_{i+j+1}^{(m)}\right)_{i,j=0}^{n-1}.$$

In general we get with induction

$$\det\left(f_{i+j+1}\right)_{i,j=0}^{n+km-1} = (-1)^{\binom{m}{2}k} a^{2kn+k^2m} \det\left(f_{i+j+1}\right)_{i,j=0}^{n-1}.$$

This implies (3.15).

**Example 3.5**

If $b_n = m+n$ for $n > 0$ then the power sequence is $m+2, m+2, 2, 2, \cdots$ and

$$f(x) = \sum_{n\geq 0} f_n x^n = \cfrac{1}{1-\cfrac{a_0x^{2+m}}{1-\cfrac{a_1x^{2+m}}{1-\cfrac{a_2x^2}{1-\cfrac{a_3x^2}{1-\ddots}}}}}.$$

Then $d(0)=1,\ d(n)=0$ for $0 < n < m+1$ and for $n > 0$

$$d(m+n) = (-1)^{\binom{m+1}{2}} a_0^{n+m} a_1^{n-1} a_2^{n-2} \cdots a_{n-1}. \tag{3.16}$$

**Example 3.6**

*Let* $(p_n)_{n\geq 0} = (1,2,2,1,2,2,1,2,2,\cdots)$. *This gives*

$(b_n)_{n\geq 0} = (0,0,2,2,3,4,5,5,7,7,8,9,10,10,12,12,13,14,\cdots)$ *with* $b_{n+6} = b_n + 5$.

*The Hankel transform is*

$1,\ 0,\ -(a_0a_1)^2,\ -(a_0a_1)^3(a_2a_3),\ -(a_0a_1)^4(a_2a_3)^2a_4,\ -(a_0a_1)^5(a_2a_3)^3a_4^2a_5,\ 0,$
$(a_0a_1)^7(a_2a_3)^5a_4^4a_5^3(a_6a_7)^2,\ (a_0a_1)^8(a_2a_3)^6a_4^5a_5^4(a_6a_7)^3(a_8a_9),\ (a_0a_1)^9(a_2a_3)^7a_4^6a_5^5(a_6a_7)^4(a_8a_9)^2a_{10},\cdots.$

We conclude this paper with some examples where the assumptions about $(b_k)$ are not satisfied.

**Example 3.7**

*Let*

$$f(x)=\cfrac{1}{1-\cfrac{x^2}{1-\cfrac{x}{1-\cfrac{x^2}{1-\ddots}}}}$$

*with power sequence* $(2,1,2,1,2,1,\cdots)$. *Then the Hankel transform is*

$$(d(n))_{n\geq 0}=(1,1,0,-1,-1-1,0,1,1,1,0,-1,\cdots). \tag{3.17}$$

**Proof**

Here we have

$f(x)=\dfrac{1}{1-x^2g(x)}$ with $g(x)=\dfrac{1}{1-xf(x)}$.

This implies $f(x)-x^2f(x)g(x)=1$ and $g(x)-xf(x)g(x)=1$ and $f(x)-1=x(g(x)-1)$.

Therefore the coefficients satisfy

$$f_{n+2}=g_{n+1}. \tag{3.18}$$

This gives by Lemma 2.1

$$\det\left(f_{i+j}\right)_{i,j=0}^{n}=\det\left(g_{i+j}\right)_{i,j=0}^{n-1}=\det\left(f_{i+j+1}\right)_{i,j=0}^{n-2}=\det\left(g_{i+j-1}\right)_{i,j=0}^{n-2}=-\det\left(f_{i+j+2}\right)_{i,j=0}^{n-4}.$$

Using (3.18) we get

$$\det\left(f_{i+j}\right)_{i,j=0}^{n}=-\det\left(f_{i+j+2}\right)_{i,j=0}^{n-4}=-\det\left(g_{i+j+1}\right)_{i,j=0}^{n-4}=-\det\left(f_{i+j}\right)_{i,j=0}^{n-4}.$$

Since the initial values are

$d(0)=1, d(1)=1, d(2)=0, d(3)=-1$ we get (3.17).

For arbitrary coefficients $a_n$ the first terms of the Hankel transform $(d(n))_{n\geq 0}$ are
$1, a_0, 0, -a_0^3a_1^2a_2^2, -a_0^4a_1^3a_2^3a_3a_4, -a_0^5a_1^4a_2^4a_3^2a_4^3, -a_0^6a_1^5a_2^5a_3^3a_4^3a_5a_6\left(a_4a_7-a_5a_6\right)$,
$a_0^7a_1^6a_2^6a_3^4a_4^4a_5^2a_6^2\left(a_4a_6^2+a_5a_6a_7a_8-a_4a_7^2a_8\right),\cdots$.
Note that in this case $(b_n)_{n\geq 0}=(0,1,1,3,2,5,3,7,4,\cdots)$.

**Example 3.8**

Finally we compute all Hankel determinants for the power sequence $\left(p_n\right)=\left(1,2,1,1,1,\cdots\right)$ which we already used for a counter example.

We write

$$f(x)=\frac{1}{1-a_0xg(x)},\quad g(x)=\frac{1}{1-a_1x^2h^{(0)}(x)},\quad h^{(0)}(x)=\frac{1}{1-a_2xh^{(1)}(x)}$$

where $g(x)$ has power sequence $\left(2,1,1,1,\cdots\right)$ and $h^{(0)}(x)$ and $h^{(1)}(x)$ have power sequence $\left(1,1,1,\cdots\right)$.

Then

$$d(n)=\det\left(f_{i+j}\right)_{i,j=0}^{n}=a_0^n\det\left(g_{i+j+1}\right)_{i,j=0}^{n-1}$$

Thus $d(0)=1$ and $d(1)=0.$

For $n>1$ we get

$$\det\left(g_{i+j+1}\right)_{i,j=0}^{n-1}=a_1^n\det\left(h_{i+j-1}^{(0)}\right)_{i,j=0}^{n-1}$$

and

$$\det\left(h_{i+j-1}^{(0)}\right)_{i,j=0}^{n-1}=-a_2^{n-2}\det\left(h_{i+j+2}^{(1)}\right)_{i,j=0}^{n-3}.$$

Therefore $d(2)=-\left(a_0a_1\right)^2$

and by (3.12) for $n\geq 0$

$$d(n+3)=-\left(a_0a_1\right)^{n+3}a_2^{n+1}a_3^{n+1}\left(a_4a_5\right)^n\left(a_6a_7\right)^{n-1}\cdots\left(a_{2n+2}a_{2n+3}\right)\sum_{i_1<i_2-1<i_3-2<\cdots<i_n-n+1\leq n}a_{i_1+3}a_{i_2+3}\cdots a_{i_n+3}. \quad (3.19)$$

If all $a_n=1$ then $h^{(0)}(x)=h^{(1)}(x)=\sum_{n\geq 0}C_nx^n.$

Therefore

$$d(n+3)=\det\left(g_{i+j+1}\right)_{i,j=0}^{n+2}=\det\left(h_{i+j-1}^{(0)}\right)_{i,j=0}^{n+2}=-\det\left(h_{i+j+2}^{(1)}\right)_{i,j=0}^{n}=-(n+2).$$

Thus in this case we get

$$\left(d(n)\right)_{n\geq 0}=\left(1,0,-1,-2,-3,\cdots\right).$$

## References

[1] George Andrews and Jet Wimp, Some $q-$orthogonal polynomials and related Hankel determinants, Rocky Mountain J. Math. 32(2), 2002, 429-442

[2] Paul Barry, On sequences with {-1,0,1} Hankel transforms, arXiv: 1205.2565

[3] Paul Barry, On the Hankel transform of C-fractions, arXiv:1212.3490

[4] Victor I. Buslaev, On Hankel determinants of functions given by their expansions in P-fractions, Ukrainian math. J. 62, 2010, 358-372

[5] Johann Cigler, Some nice Hankel determinants, arXiv:1109.1449

[6] Johann Cigler, Ramanujan's $q$ – continued fractions and Schröder-like numbers, arXiv: 1210.0372

[7] Evelyn Frank, Orthogonality properties of C-fractions, Bull. AMS 55 (1949), 384-390

[8] Christian Krattenthaler, Advanced determinant calculus, Séminaire Lotharingien Combin. 42 (1999), B42q

[9] Oskar Perron, Die Lehre von den Kettenbrüchen, B.G. Teubner 1929